\documentclass[12pt, reqno]{amsart}
\usepackage{amsmath, amsthm, amscd, amsfonts, amssymb, graphicx, color, enumerate}
\usepackage[bookmarksnumbered, colorlinks, plainpages]{hyperref}

\makeatletter \oddsidemargin.9375in \evensidemargin \oddsidemargin
\def\Speaker{$^{*}$\protect\footnotetext{$^{*}$ S\lowercase{peaker.}}}
\def\authorsaddresses#1{\dedicatory{#1}}
\newtheorem{theorem}{Theorem}[section]

\newtheorem{proposition}[theorem]{Proposition}
\newtheorem{corollary}[theorem]{Corollary}
\theoremstyle{definition}

\theoremstyle{remark}

\numberwithin{equation}{section}

\begin{document}
\setcounter{page}{1}

\noindent {\footnotesize The Extended Abstracts of \\
The 4$^{\rm th}$ Seminar on Functional Analysis and its Applications\\
2-3rd March 2016, Ferdowsi University of Mashhad, Iran}\\[1.00in]

\title[On certain product of Banach modules]{On certain product of Banach modules}

\author[M. Ramezanpour]{Mohammad Ramezanpour$^1$\Speaker}

\authorsaddresses{$^1$ School of Mathematics and Computer Science,
Damghan University, P. O. Box 36716, Damghan 41167, Iran.\\
ramezanpour@du.ac.ir\\
}


\begin{abstract}
Let $A$ and $B$ be Banach algebras and let $B$ be an algebraic Banach $A-$bimodule.  
Then the $\ell^1-$direct sum $A\times B$ equipped with   the multiplication
$$(a_1,b_1)(a_2,b_2)=(a_1a_2,a_1\cdot b_2+b_1\cdot a_2+b_1b_2),~~ (a_1, a_2\in A, b_1, b_2\in B)$$
is  a Banach algebra denoted by $A\bowtie B$. Module extension algebras, Lau product  and also the  
direct sum of  Banach algebras are the main examples satisfying this framework. 
We obtain characterizations of bounded approximate identities, spectrum,
and topological center of this product.
This provides a unified  approach for obtaining some known
results of both module extensions and  Lau product of Banach algebras.
\end{abstract}

\maketitle


\section{Introduction}

Let $A$ and $B$ be Banach algebra and $B$ is a Banach $A-$bimodule,
we say that $B$ is an algebraic Banach $A-$bimodule if
$$a\cdot(b_1b_2)=(a\cdot b_1)b_2,\quad
(b_1b_2)\cdot a=b_1(b_2\cdot a),\quad
(b_1\cdot a)b_2=b_1(a\cdot b_2),$$
for each $b_1, b_2\in B$ and $a\in A$.
Then  the Cartesian product
$A \times B$ with the algebra multiplication
$$(a_1,b_1)(a_2,b_2)=(a_1a_2,a_1\cdot b_2+b_1\cdot a_2+b_1b_2),$$
and with the norm $\|(a,b)\|=\|a\|+\|b\|$ becomes a Banach algebra
provided $\|a\cdot b\|\leq\|a\|\|b\|$, which we denote it by
 $A\bowtie B$.
We note that if  we identify $A\times\{0\}$ with $A$, and $\{0\}\times B$ with $B$, in $A\bowtie B$,
then $B$ is a closed ideal while $A$ is a closed subalgebra of $A\bowtie B$, and
$(A\bowtie B)/B$ is isometrically isomorphic to $A$.
In other words, $A\bowtie B$ is a strongly splitting Banach algebra extension of
$A$ by $B$.

Besides giving a new
method of constructing Banach algebras, the product $\bowtie$ has relevance with the following 
known products.

\begin{enumerate}[(a)]
  \item (\textit{Direct product of  two Banach algebras}) 
  Let $A$ and $B$ be Banach algebra. If we define 
  $a\cdot b=b\cdot a=0$, then $B$ is an algebraic Banach $A$-bimodule and 
  $$(a_1,b_1)(a_2,b_2)=(a_1a_2,b_1b_2).$$
  Therefore, $A\bowtie B$ is the direct product $A\times_1 B$.
  
  \item (\textit{The module extensions}) 
  Let $X$ be a Banach $A$-bimodule. Define $x_1x_2=0$, then
  $X$ is an algebraic Banach $A$-bimodule, and 
  $$(a_1,x_1)(a_2,x_2)=(a_1a_2,a_1\cdot x_2+x_1\cdot a_2).$$
  Therefore, $A\bowtie X$ is the module extension $A\oplus_1 X$. 
  Module  extensions are known as a rich source of (counter-)examples in 
  various situations in abstract harmonic analysis and functional analysis, \cite{Zhang.wameba}.
  
  \item (\textit{$\theta-$Lau  product of Banach algebras}) 
  Let $A$ and $B$ be Banach algebra and $\theta\in \Delta(A)$, the set of all  
  non zero multiplicative linear functional on $A$. Then $B$ with a module actions given by
  $a\cdot b=b\cdot a=\theta(a)b$ is an algebraic Banach $A$-bimodule and
  $$(a_1,b_1)(a_2,b_2)=(a_1a_2, \theta(a_1)b_2+\theta(a_2)b_1+b_1b_2).$$
  Thus $A\bowtie B$ is the $\theta$-Lau product $A{~}_{\theta}\!\!\times B$.
  This   product   was introduced by Lau \cite{L.laualg} for certain class of Banach algebras  
   and followed by Sangani Monfared \cite{SM.laupro} for the general case. An elementary  very 
   familiar example is  the case that  $A=\mathbb{C}$ with $\theta$ as the identity 
   character $i$ that we get the unitization $B^\sharp=\mathbb{C}{~}_{\theta\!\!}\times B$ of  $B$.
   
  \item (\textit{$T-$Lau  product of Banach algebras})
  Let $A$ and $B$ be Banach algebra and $T:A\to B$ be an algebra homomorphism with $\|T\|\leq 1$.
  Define  $a\cdot b=T(a)b=b\cdot a$. Then $B$ is an algebraic Banach $A$-bimodule   and
  $$(a_1,b_1)(a_2,b_2)=(a_1a_2,T(a_1)b_2+b_1T(a_2)+b_1b_2).$$
  Thus $A\bowtie B$ is the $T-$Lau product $A{~}_T\!\!\times B$.
  This type of product was first introduced by
   Bhatt and Dabhi for the case where $B$ is commutative and was extended by
  Javanshiri and Nemati for the general case \cite{Javn-Nemati.OcpBasp}. 
\end{enumerate}

The purpose of the present note is to determine the Gelfand space of $A\bowtie B$
which turns out to be non trivial even though $A\bowtie B$ need not be commutative
and to discuss the topological center  of $A\bowtie B$. These
topics are central to the general theory of Banach algebras.

\section{Main results}

Let $B$ be a Banach $A-$bimodule, we recall that 
$B$ is called symmetric if $a\cdot b=b\cdot a$ for all $a\in A$ and $b\in B$.

We start with the following propositions which characterize the basic properties
of $A\bowtie B$ in terms of $A$ and $B$. These results extend related results in 
\cite{Javn-Nemati.OcpBasp, SM.laupro}.
\begin{proposition}\cite{Ram}
Let $B$ be an algebraic Banach $A$-bimodule. Then 
$A\bowtie B$ is commutative if and only if $B$ is a symmetric Banach
$A$-bimodule and both $A$ and $B$ are commutative.
\end{proposition}

\begin{proposition}\cite{Ram}
Let $B$ be an algebraic Banach $A$-bimodule. Then 
 $(a_0, b_0)$ is an  identity for $A\bowtie B$
  if and only if $a_0$ is an identity for $A$, $b_0\cdot a=a\cdot b_0=0$ for all $a\in A$ and
  $a_0\cdot b+b_0b=b\cdot a_0+bb_0=b$ for all $b\in B$.
\end{proposition}

\begin{proposition}\cite{Ram}
Let $B$ be an algebraic Banach $A$-bimodule.
Then $\{(a_\alpha, b_\alpha)\}$ is a bounded left approximate identity for $A\bowtie B$
  if and only if $\{a_\alpha\}$ is a bounded left approximate identity for $A$, $\|b_\alpha\cdot a\|\to 0$
  for all $a\in A$ and $a_\alpha\cdot b+b_\alpha b\to b$ for all $b\in B$.
\end{proposition}

The dual of the space $A\bowtie B$ can be identified with $A^*\times B^*$ in the natural
way $(\varphi, \psi)(a,b)=\varphi(a)+\psi(b)$. The dual norm on $A^*\times B^*$ is of course
the maximum norm $\|(\varphi,\psi)\|=\max\{\|\varphi\|, \|\psi\|\}$.  
The following result identifies the Gelfand space of $A\bowtie B$.
This  is a generalization of \cite[Theorem 2.2]{Javn-Nemati.OcpBasp} and \cite[Proposition 2.4]{SM.laupro}.
\begin{proposition}\cite{Ram}
Let $B$ be an algebraic Banach $A$-bimodule. If $E:=\{(\varphi,0) : \varphi\in \Delta(A)\}$ and
$F:=\{ (\varphi, \psi) : \varphi\in \Delta(A)\cup\{0\},~ \psi\in \Delta(B),\text{ and} ~
a\cdot \psi=\psi\cdot a=\varphi(a)\psi~~ \forall a\in A \},$
then $\Delta(A\bowtie B)=E\cup F$.
\end{proposition}

\begin{corollary}
Let $A$ and $B$ be commutative Banach algebras and $B$ is an  algebraic Banach $A$-bimodule which is also symmetric.
Then $A\bowtie B$ is semisimple if and only if both $A$ and $B$ are semisimple.\\
\end{corollary}

Let $X$ be a Banach $A$-bimodule, for $a\in A, x\in X, x^*\in X^*, a^{**}\in A^{**}$
and $x^{**}\in X^{**}$ we define
$$\begin{array}{lcl}
(x^{**}\circ a^{**})(x^*)=x^{**}(a^{**}\circ x^*),&\quad\qquad
&(a^{**}\circ x^{**})(x^*)=a^{**}(x^{**}\circ x^*)\\
(a^{**}\circ x^*)(x)=a^{**}(x^*\circ x),&
&(x^{**}\circ x^*)(a)=x^{**}(x^*\circ a)\\
(x^*\circ x)(a)=x^*(x\cdot a),& &(x^*\circ a)(x)=x^*(a\cdot x).
\end{array}$$
Clearly, for each $a^{**}\in A^{**}$ and $x^{**}\in X^{**}$ the mappings
$b^{**}\to b^{**}\circ x^{**}:A^{**}\to X^{**}$ and
$y^{**}\to y^{**}\circ a^{**}:X^{**}\to X^{**}$ are w$^*$-w$^*$-continuous.
The first topological centres of module actions of
$A$ on $X$  may therefore be defined as
\begin{align*}
&Z^1_A(X^{**})=\{x^{**}\in X^{**} : a^{**}\to x^{**}\circ a^{**}
\mbox{  is w$^*$-w$^*$-continuous}\},\\
&Z^1_X(A^{**})=\{a^{**}\in A^{**}~:~x^{**}\to a^{**}\circ x^{**}
\mbox{  is w$^*$-w$^*$-continuous}\}.
\end{align*}

If we consider $X=A$ with the natural $A$-bimodule structure then we obtain the
first Arens product  on $A^{**}$. In this case we write $Z^1(A^{**})$ 
instead $Z^1_A(A^{**})$.
The Banach algebra $A$ is called Arens regular if $Z^1(A^{**})=A^{**}$.

To state our next result we note that if $B$ is an algebraic Banach $A$-bimodule then
$B^{**}$ is an algebraic Banach $A^{**}$-bimodule, when
$A^{**}$ and $B^{**}$  are equipped with their
first Arens products.
\begin{theorem}\cite{Ram}
Let $B$ be  an algebraic Banach $A$-bimodule.
\begin{enumerate}
  \item Suppose that $A^{**}, B^{**}$, and $(A\bowtie B)^{**}$ are equipped with their first
  Arens products. Then
  $$(A\bowtie B)^{**}\cong A^{**}\bowtie B^{**}\qquad(\text{isometric~isomorphism}).$$
  \item If $B$ is Arens regular then 
  \begin{align*}
  &Z^1((A\bowtie B)^{**})=\left(Z^1(A^{**})\cap Z^1_B(A^{**})\right)\times Z^1_A(B^{**})
  \end{align*}
\end{enumerate}
\end{theorem}

%
%

\bibliographystyle{amsplain}

\end{document}